\documentclass[12pt]{amsproc}
\usepackage{amssymb}

\usepackage{amsfonts,amsmath,amsthm}
\usepackage{amssymb}
\usepackage[utf8]{inputenc}
\usepackage[english]{babel}
\usepackage{graphicx}
\usepackage{epstopdf}
\usepackage{tikz}
\usetikzlibrary{shapes,arrows,backgrounds}
\usepackage{subcaption}
\usepackage{array}
\usepackage{url}
\usepackage{booktabs}

\newtheorem{theorem}{Theorem}[section]
\newtheorem*{theorem-non}{Theorem}
\newtheorem{lemma}[theorem]{Lemma}

\newtheorem{proposition}[theorem]{Proposition}
\newtheorem*{proposition-non}{Proposition}
\newtheorem{prop}[theorem]{Proposition}

\newtheorem*{conjecture-non}{Conjecture}
\theoremstyle{remark}
\newtheorem{nota}[theorem]{Remark}

\newtheorem{remark}[theorem]{Remark}

\def\equi{=}
\def\interior{\operatorname{int}}
\def\conv{\operatorname{conv}}

\def\R{\mathbb R}
\def\Z{\mathbb Z}
\def\N{\mathbb N}
\def\B{\mathcal B}

\def\norma{\|\cdot\|}

\title{A characterisation of Euclidean normed planes via bisectors}

\begin{document}

\author
{Javier~Cabello~Sánchez, Adri\'an~Gordillo-Merino}
\address{Departamento de Matem\'{a}ticas, Universidad de Extremadura, Avda. de Elvas s/n, 06006 Badajoz. Spain}
\email{coco@unex.es, adgormer@unex.es}
\thanks{First author supported in part by DGICYT project MTM2016$\cdot$76958$\cdot$C2$\cdot$1$\cdot$P (Spain) and Junta de Extremadura programs GR$\cdot$15152 and  IB$\cdot$16056.} 
\thanks{Second author has been partially supported by Junta de Extremadura and FEDER funds.}


\begin{abstract}
Our main result states that whenever we have a non-Euclidean norm
$\norma$ on a two-dimensional vector space $X$, there exists some $x\neq 0$ such that for every $\lambda\neq 1, 
\lambda>0$, there exist $y, z\in X$ verifying that $\|y\|=\lambda\|x\|$, $z\neq 0$, and $z$ belongs to the bisectors
$B(-x,x)$ and $B(-y,y)$. Throughout this paper we also state and prove some other
simple but maybe useful results about the geometry of the unit sphere of strictly convex planes.
\end{abstract}

\maketitle

\section{Introduction}

\bigskip

In a normed linear space $(X,\norma)$, a vector $x$ is said to be isosceles orthogonal to 
a vector $y$ (denoted by $x\perp_I y\,$) if $\|x-y\|=\|x+y\|$. Isosceles orthogonality was 
introduced by R.C.~James in \cite{james1945}. Since then, several papers and surveys have 
studied properties related to the geometric structure of the space in the light of that 
notion of orthogonality, and various characterisation results (e.g., for strict convexity) 
have been obtained. Two interesting surveys on this topic are \cite{alonso2012birkhoff} and \cite{MARTINI200493}, and the monograph \cite{thomson1996} gives further background.

\medskip

In this paper, $(X,\norma)$ will denote a 2-dimensional normed space (usually referred to as 
a {\it Minkowski plane}), and $S_X\,$ and $B_X\,$ will stand for the {\it unit sphere} and the {\it closed 
unit ball}, respectively. We hasten to remark that, as we are dealing with normed spaces, $B_X\,$ is always a planar convex body centred at the origin and $S_X\,$ agrees with its boundary. The segment 
joining two points $x, y$ will be denoted as $[x, y]$. As we will deal with segments, 
intervals and two-dimensional vectors, we need to determine the meaning of $(x,y)$. 
Throughout the paper, this will denote a vector in $X$. Of course, we will need to have a basis $\{e,v\}$ fixed previously, so that $(x,y)$ means $xe+yv$. For open intervals 
(or segments) we will use the notation $]x,y[$, and for semiopen intervals 
(segments) we will write $[x,y[$ and $]x,y]$. For the linear span of a couple of vectors $x,y\in X$ 
we will use $\langle x,y\rangle$. 

We will utilize the concept of bisector of the segment joining two points. For $x,y\in X\,$, 
the bisector of $[x,y]\,$ is defined in this way 
(check \cite{alonso2012birkhoff, barschnei2014, BisectorsHMW, BisectorsJS, ma2000bisectors} for instance):
$$
B(x,y)=\{z\in X:\|x-z\|=\|y-z\|\}\,.
$$

\medskip

In Section \ref{sideresults} we prove Proposition \ref{anticonjetura}, now stated in a slightly different way:

\begin{proposition-non} A norm $\norma$ on $X\,$ is strictly 
convex if and only if for every nonzero $z\in X\,$ there exists, up to $\pm 1$, exactly one vector which is isosceles orthogonal to $z$ in $S_X$.
\end{proposition-non}

This solves in the negative the following conjecture, proposed by Alonso, Martini and Wu (\cite[Conj. 5.3.]{alonso2012birkhoff}),  
with a different approach to the one used to give the solution that can be found in \cite[Prop. 5]{alonso1994} combined with \cite[Cor. 2.5]{jiliwu2011}:

\begin{conjecture-non}
In any non-Euclidean Minkowski plane X, there exist $x, y\in S_X$, 
with $x\neq  \pm y$, such that $B(-x,x)\cap B(-y,y)\neq\{0\}.$ 
\end{conjecture-non}

As far as we know, this if-and-only-if statement of our Proposition \ref{anticonjetura} cannot be found in the literature.

\medskip

Section \ref{mainresult} is mainly devoted to the main result (Theorem \ref{maintheorem}) in our paper. We propose this characterisation of Euclidean normed planes:

\begin{theorem-non}
The norm $\norma\,$ is not Euclidean if and only if, for some $x\neq 0\,$ and each 
$\lambda\in (0,+\infty)\setminus \{1\}\,$, there exists $y\,$ such that 
$\|y\|=\lambda \|x\|\,$, $\langle x,y\rangle =X\,$, and $B(-x,x)\cap B(-y,y)\neq 0\,$.
\end{theorem-non}

\bigskip
 
\begin{nota}\label{lalinealidaddelbisector}
Please observe that the definitions of isosceles orthogonality and bisectors recalled above lead to the following equivalences: 

\noindent $z\in B(-x,x)$ if and only if $x\perp_I z$ if and only if $x\in B(-z,z)$. 

Besides, it is easily checked that bisectors enjoy a certain property of linearity:
$$
B(\lambda x+z,\lambda y+z)=z+\lambda B(x,y)\,,\,\,\,\forall x,y,z\in X\,,\,\,\,\forall \lambda \in \mathbb{R}\,.
$$

Therefore, for any $a,b\in X\,$, 
$$
B(a,b)=\frac{a+b}{2}+\frac{\|a-b\|}{2}B\left(\frac{a-b}{\|a-b\|},\frac{b-a}{\|b-a\|}\right)\,,
$$
so the geometric properties of bisectors can be determined just by careful analysis of properties of bisectors of the type $B(-x,x)\,$, with $x\in S_X\,$.
\end{nota}

\medskip

\section{Side results}\label{sideresults}

We will prove that the conjecture in \cite[Conj. 5.3.]{alonso2012birkhoff} 
is false by proving Proposition \ref{anticonjetura}. We hasten to remark 
that one implication can be seen in~\cite[Proposition 5]{alonso1994}, 
while the other is proven in \cite[Corollary 2.5]{jiliwu2011}. This is going to be 
achieved in a different, say more geometric, manner to those. 
In order to complete the proof, some basic results are presented.

\begin{proposition}\label{anticonjetura}
Let $(X,\|\cdot\|)$ be a normed plane. Then, it is strictly convex if and only if 
$B(-x,x)\cap B(-y,y)=0$ for every linearly independent $x, y\in S_X$. 
\end{proposition}

\begin{nota}
As we have noted in Remark~\ref{lalinealidaddelbisector}, every bisector is an affine 
transformation of a $B(-x,x)$ for some $x\in S_X$. This readily implies that the following 
conditions are equivalent:
\begin{itemize}
\item There exists $\lambda>0$ such that $B(-x,x)\cap B(-y,y)=0$ for every linearly independent $x, y\in\lambda S_X$.
\item $B(-x,x)\cap B(-y,y)=0$ for every linearly independent $x, y\in S_X$.
\item For every $\lambda>0$ and every linearly independent $x, y\in\lambda S_X$ one has $B(-x,x)\cap B(-y,y)=0$.
\item For every $\lambda>0, z\in X$ and every linearly independent $x, y\in\lambda S_X$ one has $B(z-x,z+x)\cap B(z-y,z+y)=z\,$.
\item $B(x,x')\cap B(y,y')=(x+x')/2\,$, whenever $\|x-x'\|=\|y-y'\|$ and $x+x'=y+y'\,$. 
\end{itemize}
\end{nota}

\begin{remark}
Our problem is to determine what happens when $0\neq z\in X, x, y\in S_X$ 
are such that $z \in B(-x,x)\cap B(-y,y)$ and $x$ and $y$ 
are linearly independent. This is equivalent to $\|z-x\|=\|z+x\|$ and $\|z-y\|=\|z+y\|$ or, 
with $\lambda_x=\|z-x\|^{-1}, \lambda_y=\|z-y\|^{-1}$, to 
$$\lambda_x(z+x), \lambda_x(z-x),\lambda_y(z+y),\lambda_y(z-y)\in S_X.$$
We have, then, two couples of points in the sphere, say 
$$a=\lambda_x(z+x), a'=\lambda_x(z-x), b=\lambda_y(z+y), b'=\lambda_y(z-y)$$ 
and two positive values (not necessarily different) $\alpha, \beta$, such that 
$$\alpha z=a+a'=\lambda_x(z+x+z-x)=2\lambda_x z, \beta z=b+b'=\lambda_y(z+y+z-y)=2\lambda_y z.$$ 
In particular, $\alpha=2\lambda_x=\|a+a'\|/\|z\|, \beta=2\lambda_y=\|b+b'\|/\|z\|$. 
In the statement we ask $x$ and $y$ to be linearly independent, and we also have 
$$a-a'=\lambda_x(z+x-(z-x))=2\lambda_xx, b-b'=\lambda_y(z+y-(z-y))=2\lambda_yy, $$ 
so $a-a'$ and $b-b'$ must be independent, too. 
From these last equalities we obtain $2\lambda_x=\|a-a'\|, 2\lambda_y=\|b-b'\|$ so mixing 
with the previous ones, we obtain $\|a-a'\|=\|a+a'\|/\|z\|, \|b-b'\|=\|b+b'\|/\|z\|$, or 
$$\frac {\|a+a'\|}{\|a-a'\|}=\|z\|=\frac {\|b+b'\|}{\|b-b'\|}.$$ 
\end{remark}

\begin{proof}[Proof of the easy implication of Prop.~\ref{anticonjetura}] 

Suppose $\norma$ not to be strictly convex. Then, there is some segment $[c,c']\subset S_X$. 
Take $$a=\frac 14(3c+c'), a'=\frac 14(-3c'-c), b=c, b'=-\frac 12(c+c').$$ It is 
straightforward that $a, a', b, b'\in S_X$, that 
$0\neq z=\frac 12(b+b')=\frac 12(a+a')$, and also that $a-a'$ and $b-b'$ are independent. 
Also, $\|a-a'\|=\|b-b'\|=2$, and this implies that 
$\frac{\|a+a'\|}{\|a-a'\|}=\frac {\|b+b'\|}{\|b-b'\|}$. 
Please observe that we have $0\neq z\in B(-x,x)\cap B(-y,y)$, with 
$x=a-z, y=b-z$. 
\end{proof}

\subsection{Proof of the other implication}

For the remainder of this section, let us suppose that $\|\cdot\|$ is strictly convex. 

\medskip

We will also assume that 
$0\neq z\in X, a, a', b, b'\in S_X$, $\beta\geq\alpha>0\,$ are such that 
$$a+a'=\alpha z\,,\, b+b'=\beta z\,,\, \frac{\|a+a'\|}{\|a-a'\|}=\frac{\|b+b'\|}{\|b-b'\|}, $$
and $a-a'$ and $b-b'$ are linearly independent. 

We will split the proof in some elementary results that may be useful for broader purposes. 

First, let us fix some notations: 

We will consider $X$ endowed with the basis $\{z, (a-a')/2\}$, so we have 
$a\equi(\alpha,1), a'\equi(\alpha,-1)$, $\delta=2/\|a-a'\|$ and 
$d\equi(0,\delta), d'\equi(0,-\delta)$, with $a, a', d, d'\in S_X$. Of course, 
$\delta>1$. Consider the lines $r^+$ and $r^-$ 
defined, respectively, by $(0,\delta), (\alpha,1)$ and $(0,-\delta), (\alpha,-1)$. They are given by $r^+(x)=\delta+(1-\delta)x/\alpha$, $r^-(x)=-\delta+(\delta - 1)x/\alpha\,$, and 
the only point that their graphs have in common is $c\equi(\delta\alpha/(\delta-1),0)$.
So, we have 
$$d\equi(0,\delta)\equi(0,r^+(0)), d'\equi(0,-\delta)\equi(0,r^-(0)),$$ 
$$a\equi(\alpha,1)\equi(\alpha,r^+(\alpha)), a'\equi(\alpha,-1)\equi(\alpha,r^-(\alpha)),$$ 
$$c\equi(\delta\alpha/(\delta-1),0)\equi(\delta\alpha/(\delta-1),r^+(\delta\alpha/(\delta-1))). $$

\begin{lemma}\label{traptrian}
Consider the convex hull $\conv\{d,d',c\}$ and the vertical line $\{\alpha\}\times\R$. The following {\em  symmetric} inclusions are verified: 
$$\conv\{d, d', c\}\cap(]0,\alpha[\times\R)\subset\interior(B_X).$$
$$B_X\cap(]\alpha,\infty[\times\R)\subset \interior(\conv\{d, d',c\}). $$
\end{lemma}

\begin{proof}
Observe that 
$$\conv\{d, d', c\}\cap(]0,\alpha[\times\R)=\conv\{d, d',a ,a'\}\cap(]0,\alpha[\times\R)\mathrm{\ and}$$ 
$$\conv\{d, d',c\}\cap(]\alpha,\infty[\times\R)=\conv\{a, a',c\}\cap(]\alpha,\infty[\times\R). $$
For the first part, take $x\equi(x_1,x_2)\in\conv\{a, a', d, d'\}$, with $x_1\in]0, \alpha[$. 
We will show that $x\in\interior(B_X)$. As 
$$x\in\conv\{(0,r^+(0)), (0,r^-(0)), (\alpha, r^+(\alpha)), (\alpha, r^-(\alpha))\},$$
we have $r^-(x_1)\leq x_2\leq r^+(x_1)$. As $\|\cdot\|$ is 
strictly convex, both $(x_1,r^+(x_1))$ and $(x_1,r^-(x_1))$ belong to the interior of $B_X$. 
As $(x_1, x_2)$ is a convex combination of $(x_1,r^+(x_1))$ and $(x_1,r^-(x_1))$, we have 
$x\in\interior(B_X)$, too. 

As for the second part, let $x\equi(x_1,x_2)\in B_X$, with $x_1>\alpha$. Suppose 
that $x_2\geq r^+(x_1)=\delta+(1-\delta)x_1/\alpha$. Then 
the strict convexity of $\|\cdot\|\,$, $\|x\|\leq 1\,$ and $\|d\|= 1\,$ imply that
$$1>\left\|\frac{\alpha}{x_1}x+\frac{x_1-\alpha}{x_1}d\right\|=
\left\|\frac{\alpha}{x_1}(x_1,x_2)+\frac{x_1-\alpha}{x_1}(0,\delta)\right\|=$$
$$=\left\|\left(\alpha,\frac{\alpha x_2}{x_1}\right)+ \left(0,\delta-\frac{\alpha\delta}{x_1}\right)\right\|
=\left\|\left(\alpha,\delta+\frac{x_2-\delta}{x_1}\alpha\right)\right\|.$$
Now, we have 
$\delta+\frac{\alpha(x_2-\delta)}{x_1}\geq \delta+\frac{\alpha}{x_1}(\delta+(1-\delta)\frac{x_1}{\alpha}-\delta)
=1$, which implies that  $(\alpha,1)$ is a convex combination of $(\alpha,-1)$ and 
$(\alpha,\delta+(x_2-\delta)\alpha/x_1)$, but $\|(\alpha,1)\|=\|(\alpha,-1)\|=1$ and 
$\|(\alpha,\delta+(x_2-\delta)\alpha/x_1)\|<1,$ a contradiction. 
The case $x_2\leq r^-(x_1)$ is analogous. 
\end{proof}

\begin{prop}\label{bolabola}
For the ball $(\alpha,0)+\delta^{-1}B_X$ we have essentially the same symmetric inclusions: 
$$B_X\cap(]\alpha, \infty[\times\R)\subset(\alpha,0)+\delta^{-1}\interior(B_X). $$
$$((\alpha,0)+\delta^{-1}B_X)\cap(]-\infty,\alpha[\times\R)\subset\interior(B_X). $$
\end{prop}

\begin{proof}
Let $x\equi(x_1,x_2)\in B_X$ be such that $x_1>\alpha$, we may suppose $x_2\geq 0$. 
Instead of showing that $x$ belongs to the interior of $(\alpha,0)+\delta^{-1}B_X$, we shall see that 
$$(\delta x_1-\delta\alpha,\delta x_2)\in \interior(B_X).$$ 
As $(x_1,x_2), (0,x_2)$ and $(0,\delta)$ belong to $B_X$, it suffices to show that 
$$(\delta x_1-\delta\alpha,\delta x_2)\in\conv\{(x_1,x_2),(0,x_2),(0,\delta)\}.$$
For this, we need $\delta x_1-\delta\alpha\in ]0,x_1[$. This is equivalent 
to $x_1<\alpha\delta/(\delta-1)$, and this inequality is true 
since $c\equi(\alpha\delta/(\delta-1),0)$ is the only point in $r^+\cap r^-$. 

As $\delta x_2>x_2$, the only we still need to show is that $(\delta x_1-\delta\alpha,\delta x_2)$ 
lies below the line defined by $(0,\delta)$ and $(x_1,x_2)$. This line is the 
graph of the function $y(t)=x_2t/x_1+\delta-\delta t/x_1$, and so we need 
$$\delta x_2<x_2(\delta x_1-\delta\alpha)/x_1+\delta-\delta (\delta x_1-\delta\alpha)/x_1.$$
After some elementary computations, we get that this inequality is equivalent to 
$$0<x_1-x_2\alpha+\alpha\delta-\delta x_1.$$

To finish the proof of the first part we only need to observe that the second 
part of Lemma~\ref{traptrian} implies that $(\alpha,1)$ is above the line 
defined by $(0,\delta)$ and $(x_1,x_2)$, so $1>x_2\alpha/x_1+\delta-\delta \alpha/x_1$. It is easy 
to see that this is also equivalent to $0<x_1-x_2\alpha+\alpha\delta-\delta x_1,$ and so we are done. 

For the second item, take $y\equi(y_1,y_2)\in(\alpha,0)+\delta^{-1}S_X$, with 
$y_1<\alpha$, and $y'\equi(y_1',y_2')\equi(2\alpha-y_1,-y_2)$ symmetric 
to $y$ with respect to $(\alpha,0)$, and suppose that $y\in S_X$. 
As both $(\alpha,1)$ and $(\alpha,0)+\delta^{-1}(\alpha,1)$ belong to 
$((\alpha,0)+\delta^{-1}S_X)\cap r^+$, there are no more points in this intersection, 
and this means that $(\alpha,0)+\delta^{-1}S_X$ lies below $r^+$ outside the 
interval $[\alpha,\alpha(1+\delta)]$. As $S_X$ lies above this line in $[0,\alpha]$, 
we get $y_1<0$. Now, $y'\in (\alpha,0)+\delta^{-1}S_X$ and $y'_1=2\alpha-y_1>2\alpha$ 
together imply $|y_2|=|y_2'|<\delta^{-1}<1$, and from this we get $y_1<-\alpha$. 

We have also $\|(y_1-\alpha,y_2)\|=\delta^{-1}$, and so we get 
$$(y_1-\alpha,y_2), (\alpha, y_2)\in \interior(B_X), (y_1,y_2)\in S_X\mathrm{\  and\ }
(y_1,y_2)\in[(y_1-\alpha,y_2), (\alpha, y_2)],$$ a contradiction. 
 
\end{proof}

\begin{lemma}\label{normamedios}
With the previous notations, $\beta>\alpha$ implies $\|b-b'\|< \|a-a'\|$. 
\end{lemma}

\begin{proof}
We may suppose that $\|z\|=1$. Recall that, in the basis we are dealing with, 
$\frac 12(a+a')\equi(\alpha,0)$ and $\frac 12(b+b')\equi(\beta,0)$. 

Let $b\equi(b_1,b_2), b'\equi(b_1',b_2')$ be the expressions in coordinates of $b, b'$ 
in the basis $\{z, (a-a')/2\}$. It is clear that $(b_1',b_2')\equi(2\beta-b_1, -b_2)$, and 
we may suppose $b_1\geq \beta>\alpha$. Then, with $\beta'=r^+(\beta)<r^+(\alpha)=1$, Lemma~\ref{traptrian} implies 
$$b\in B_X\cap([\beta,\infty[\times\R)\subset [(\beta,0)+\beta'B_X]\cap([\beta,\infty[\times\R),$$
so $\left\|b-(\beta,0)\right\|<\|a-(\alpha,0)\|$, and we are done. 
\end{proof}

\begin{lemma}\label{proyeccion}
Let $z\in \interior(B_X)\setminus\{0\}$. There exists exactly one couple $x, x'\in S_X$ such that 
$z=\frac 12(x+x')$. 
\end{lemma}

\begin{proof}
For the existence, we will define some auxiliary functions. For $t\in[0,2\pi]$, 
let $x(t)$ be defined as the only point in $S_X\cap\{\lambda(\cos(t),\sin(t)):\lambda\in]0,\infty[\}$. 
Take $z\in \interior(B_X)\setminus\{0\}$ and $f(t)$ be defined as $\|z-y(t)\|$, where 
$y(t)$ is the only point in $S_X\cap\{z+\lambda x(t):\lambda\in]0,\infty[\}$. 
It is pretty clear that all these functions are continuous and, moreover, 
$f(2\pi)=f(0)$. So, there exists $t\in[0,\pi[$ such that $f(t+\pi)=f(t)$. 
For this $t$, we have $z=\frac 12(y(t)+y(t+\pi))$. 

For the uniqueness, suppose that we have four different points 
$x, x', y, y'\in S_X$ such that $x+x'=y+y'=2z$, and take as a basis 
$\{z,1/2(x-x')\}$, so that $x+x'=y+y'\equi(2,0)$, $x\equi(1,1), x'\equi(1, -1)$ and 
$y\equi(y_1,y_2), y'\equi(y'_1y'_2)$. As usual, $\delta=1/\|(0,1)\|.$

Now, suppose $y_1>1$. By the first part of Proposition~\ref{bolabola}, $y\in S_X$ 
implies $y\in(1,0)+\delta^{-1}\interior(B_X)$. But the second part of the same Lemma 
implies that, then, $y'\in\interior(B_X)$, so we are done. 
\end{proof}

To finish our proof of the remaining implication of Proposition \ref{anticonjetura}, we only need to notice that, for $\beta>\alpha>0\,$, Lemma \ref{normamedios} leads us to a contradiction with our initial assumptions and, in case $\alpha=\beta\,$, the contradiction springs up from Lemma \ref{proyeccion}.

\medskip

\section{Main result}\label{mainresult}

Now we can state and prove the last step before the main result. 

Please recall that we are not assuming $(X,\norma)$ to be strictly convex anymore. 

\begin{prop}\label{DosRectas}
Let $x, y, z$ be nonzero vectors in $X$ and $(\gamma_n), (\delta_n)\subset\R$ be a 
couple of positive sequences converging monotonically to 0. If $\gamma_nx, 
\delta_ny\in B(-z,z)$ for every $n\in\N$, then $y=\pm x$. 
\end{prop}

\begin{proof}
We may suppose $\|z\|=1$. 

If the result does not hold, then we may take $\{x,y\}$ as a basis of $X$ and, 
in coordinates, we have $x\equi(1,0), y\equi(0,1), z\equi(z_1,z_2)$, and we may suppose $z_1, z_2>0$. 
Indeed, if $z_1<0$ then we may take $-x$ instead of $x$ and the case $z_1=0$ is absurd. 

As $\gamma_nx, \delta_ny\in B(-z,z)$, in coordinates we have 
$$\|(z_1+\gamma_n,z_2)\|=\|(z_1-\gamma_n,z_2)\|\qquad \|(z_1,z_2+\delta_n)\|=\|(z_1,z_2-\delta_n)\|, \forall\,n.$$

It will be better to have a shorter way to refer to these values, so let $\alpha_n=\|(z_1+\gamma_n,z_2)\|^{-1}$ 
and $\beta_n=\|(z_1,z_2+\delta_n)\|^{-1}$. 

Observe that we have $\alpha_n\to 1, \beta_n\to\|z\|= 1$, and also that the convexity 
of $\norma$ implies that $\|( z_1\pm\gamma_n,z_2)\|$ and $\|( z_1,z_2\pm\delta_n)\|$ 
are not lower than 1, so $\gamma_n\leq 1,\ \delta_n\leq 1$ for every $n$.

The choice of $(\alpha_n)$ and $(\beta_n)$ gives 
$$\alpha_n( z_1\pm\gamma_n,z_2), \beta_n( z_1,z_2\pm\delta_n)\in S_X, \forall\,n,$$
so we have a couple of sequences $(\alpha_n)z$ and $(\beta_n)z$ that converge to $z$ and such that
each $\alpha_nz$ is the midpoint of the segment $(\{\alpha_n z_1\}\times\R)\cap B_X$ and $\beta_nz$ 
is the midpoint of the segment $(\R\times\{\beta_nz_2\})\cap B_X$. 

Now, we are going to analyse the shape of the unit ball $B_X$ of such a norm and will 
eventually rule out every possibility. 

Suppose there is some $v\equi(v_1,v_2)\in B_X$ with $v_1>z_1$. Then, $B_X$ contains the 
triangle $\conv\{(v_1,v_2), (0,0),(z_1,z_2)\}$. On the one hand, it is actually a triangle 
provided $v_1z_2\neq z_1v_2$, but $v_1z_2= z_1v_2$ is absurd. 
On the other hand, if $v_1z_2>v_2z_1$ (respectively $v_1z_2<v_2z_1$), then the interior of 
this triangle contains $(z_1+\gamma_n,z_2)$ (resp. $(z_1-\gamma_n,z_2)$) for infinitely many $n$. 
As $( z_1\pm\gamma_n,z_2), ( z_1,z_2\pm\delta_n)\not\in\interior B_X$ for every $n$, 
and applying the analogous reasoning to $v_2$ we get $v_1\leq z_1$ and $v_2\leq z_2$ 
for every $(v_1,v_2)\in B_X$. So, we have
$$\{(v_1,v_2)\in B_X:v_1, v_2\geq 0\}\subseteq\conv\{(z_1,z_2),(0,z_2),(0,0),(z_1,0)\}.$$

Consider now $r_2$ as the line that is {\em vertically symmetric} to 
$r_1=\{(t,z_2):t\in\R\}$ with respect to $r_0=\{(t,tz_2/z_1):t\in\R\}=\langle z\rangle$, i.e., 
the line $r_2=\{(t,2tz_2/z_1-z_2)\}$: for every $t$, the (unique) point in 
$r_2\cap(\{t\}\times\R)$ is symmetric to (the point in) $r_1\cap(\{t\}\times\R)$ 
with respect to $r_0\cap(\{t\}\times\R)$. As $B_X$ lies below $r_1$, it is readily seen that 
$\beta_n(z_1,z_2-\delta_n)$ lies above $r_2$ for every $n$. So, with an argument similar to that in 
the previous paragraph, we can see that every $(v_1,v_2)\in B_X$ lies above $r_2$. 
As $r_2$ contains both $(z_1,z_2)$ and $(0,-z_2)$, the point where $r_2$ intersects 
the horizontal axis is $(z_1/2,0)$, so we can describe the situation this way:
$$\{(v_1,v_2)\in B_X:v_1, v_2\geq 0\}\subseteq\conv\{(z_1,z_2),(0,z_2),(0,0),(z_1/2,0)\}.$$

Now, consider $r_3$ as the line that is {\em horizontally symmetric} to $r_2$ 
with respect to $r_0$, i.e., the midpoint of $r_3\cap(\R\times \{t\})$ and $r_2\cap(\R\times \{t\})$ 
is $r_0\cap(\R\times \{t\})$ for every $t$. The same argument implies that $B_X$ lies below 
$r_3$, and $(0,z_2/3)$ is the point in which $r_3$ intersects the vertical axis. So, 
$$\{(v_1,v_2)\in B_X:v_1, v_2\geq 0\}\subseteq\conv\{(z_1,z_2),(0,z_2/3),(0,0),(z_1/2,0)\}.$$

If we keep iterating the process, for $n$ we get 
$$\{(v_1,v_2)\in B_X:v_1, v_2\geq 0\}\subseteq\conv\{(z_1,z_2),(0,z_2/(2n-1)),(0,0),(z_1/2n,0)\}$$
and this is absurd. So, we have finished the proof. 
\end{proof}

Finally, we have this new characterisation (in the negative) of the Euclidean case among all Minkowski planes. 
We thought about stating it in the positive, but it seemed less intelligible to us. 

\begin{theorem}\label{maintheorem}
The norm $\norma\,$ is not Euclidean if and only if, for some $x\neq 0\,$ and each $\lambda\in (0,+\infty)\setminus \{1\}\,$, there exists $y\,$ such that 
$\|y\|=\lambda \|x\|\,$, $\langle x,y\rangle =\mathbb{R}^2\,$, and $B(-x,x)\cap B(-y,y)\neq 0\,$.
\end{theorem}

\begin{proof}

$(\Leftarrow)\,$: This is the simple part: if $\norma\,$ is Euclidean, then $B(-x,x)\,$ and $B(-y,y)\,$ are both straight lines, and they are different provided that $\langle x,y\rangle=\mathbb{R}^2\,$.

$(\Rightarrow)\,$: First, let us assume $\norma\,$ not to be strictly convex. 

Let $a, b\in S_X$ such that $[a,b]\subset S_X$. Take $x=(3a+b)/4$ and $z=(a+3b)/4$ and observe 
that $(x-z)/2=(a-b)/4$ belongs both to $B(-x,x)$ and $B(-z,z)$. Furthermore, for every $\alpha\in[-1,1]$ we have 
$$\|z-\alpha(x-z)/2\|=\|(a+3b)/4-\alpha(a-b)/4\|=\|(1-\alpha)a/4+(3+\alpha)b/4\|=1$$
because the last is a convex combination of $a$ and $b$. So, the full segment $[(z-x)/2,(x-z)/2]$ 
lies in $B(-z,z)$. By symmetry, it is included in $B(-x,x)$, too. 

Now, let $\lambda>0$. If $\lambda\leq 1$, then $\lambda(z-x)/2\in B(-x,x)\cap B(-\lambda z,\lambda z)$. 
If $\lambda\geq 1$, then $(z-x)/2\in B(-x,x)\cap B(-\lambda z,\lambda z)$, so in any case the result follows with 
$y=\lambda z$.

\medskip

As for the case in which $\norma\,$ is strictly convex, suppose that for every $x\neq 0\,$ and a certain $\lambda\in (0,+\infty)\setminus \{1\}\,$, there exists no $y\,$, $\|y\|=\lambda\|x\|\,$, verifying $\langle x,y\rangle =\mathbb{R}^2\,$ and $B(-x,x)\cap B(-y,y)\neq 0\,$.

Consider $x_0\,$ such that $B(-x_0,x_0)\,$ is not a straight line (the existence of such bisector is guaranteed in any non-Euclidean Minkowski plane, see \cite[Th. 5.5.]{alonso2012birkhoff}). Then, $B(-x_0,x_0)=B(-\lambda x_0,\lambda x_0)\,$: 

Let us assume this is not the case: let $p\in B(-x_0,x_0)\,$, $p\notin B(-\lambda x_0,\lambda x_0)\,$. Now, if we take the unique (see \cite[Cor. 2.5]{jiliwu2011}) $y=B(-p,p)\cap \lambda \|x_0\|S_{X}\,$, we come to a contradiction, as $p\in B(-x_0,x_0)\cap B(-y,y)\,$.

On the other hand, $B(-\lambda x_0,\lambda x_0)=\lambda B(-x_0,x_0)\,$, as we said in Remark~\ref{lalinealidaddelbisector}; therefore, $B(-x_0,x_0)=\lambda B(-x_0,x_0)\,$.
Now, take linearly independent $e, v\in B(-x_0,x_0)$. As $B(-x_0,x_0)=\lambda B(-x_0,x_0),$ 
we have $\lambda u\in B(-x_0,x_0)$ for every $u\in B(-x_0,x_0)$ so $\lambda^n u\in B(-x_0,x_0)$ and, 
in particular, 
\begin{equation}\label{lambdaestan}
\lambda^n e, \lambda^n v\in B(-x_0,x_0)\mathrm{\ for\ every\ } n\in\Z.
\end{equation}
As $\B=\{e,v\}$ is a basis in $X$, we may take coordinates and we get 
$e\equi(1,0), v\equi(0,1), x_0\equi(\alpha,\beta)$. For the sake of clarity, we will suppose 
$x_1,x_2>0$ ---if we had $x_1<0$ we could just take $-e$ instead of $e$ and the case 
$x_1=0$ is absurd. 
We may also suppose $\lambda\in]0,1[$. Indeed, if we have $\lambda>1$ 
we can take $\mu=\lambda^{-1}$ and we are exactly in the same situation as before: 
$B(-x_0,x_0)=\mu B(-x_0,x_0)$ with $\mu<1$. 

Rewriting~(\ref{lambdaestan}) in coordinates we get: 
$$
\|(\alpha-\lambda^n,\beta)\|=\|(\alpha+\lambda^n,\beta)\|=\|(-\alpha+\lambda^n,-\beta)\|=\|(-\alpha-\lambda^n,-\beta)\|, 
$$
$$
\|(\alpha,\beta-\lambda^n)\|=\|(\alpha,\beta+\lambda^n)\|=\|(-\alpha,-\beta+\lambda^n)\|=\|(-\alpha,-\beta-\lambda^n)\|, 
$$

By Proposition \ref{DosRectas}, this cannot happen and we have finished the proof. 

\end{proof}

%

\end{document}